\documentclass[nopreprintline,12pt,times]{elsarticle}

\usepackage{amssymb}
\usepackage{amsmath}

\usepackage{bbm}
\usepackage{lineno}
\usepackage{lipsum}

\journal{Applied Numerical Mathematics}

\begin{document}
	
	\begin{frontmatter}
		
		
		
		\title{A DEEP LEARNING APPROXIMATION OF NON-STATIONARY
			SOLUTIONS TO WAVE KINETIC EQUATIONS\tnoteref{t1}}
		\tnotetext[t1]{S.W. was supported by  U.S. National Science Foundation under the SMU Research Training Group grant DMS-1840260. M.-B. T. was  funded in part by  the  NSF Grant DMS-1854453,    Humboldt Fellowship,   NSF CAREER  DMS-2044626, NSF Grant DMS-2204795. A.B. was  funded in part by  the  NSF Grant DMS-1905449,    NSF Grant DMS-2204795 and a grant from the SAR Hong Kong RGC GRF 14301321.}
		
		\author[1]{Steven Walton\corref{cor1}%
		}
		\ead{stevenw@smu.edu}

		\author[3]{Minh-Binh Tran}
		\ead{minhbinh@tamu.edu}

		\author[4,5]{Alain Bensoussan}
		\ead{alain.bensoussan@utdallas.edu, abensous@cityu.edu.hk}
		
		\cortext[cor1]{Corresponding author}
		
		\affiliation[1]{organization={Southern Methodist University},
			city={Dallas},
			state={Texas},
			country={USA}}
		\affiliation[3]{organization={Texas A\&M University},
			city={College Station},
			state={Texas},
			country={USA}}
		\affiliation[4]{organization={University of Texas Dallas},
			city={Dallas},
			state={Texas},
			country={USA}}
		\affiliation[5]{organization={City University of Hongkong},
			country={Hongkong}}

	
	
	
	\begin{abstract}
		We present a deep learning approximation, stochastic optimization based, method for wave kinetic equations.  To build confidence in our approach, we apply the method to a Smoluchowski coagulation equation with multiplicative kernel for which an analytic solution exists.  Our deep learning approach is then used to approximate the non-stationary solution to a 3-wave kinetic equation corresponding to acoustic wave systems. To validate the neural network approximation, we compare the decay rate of the total energy with previously obtained theoretical results. A finite volume solution is presented and compared with the present method.     
		
	\end{abstract}
	
	
	
	
\end{frontmatter}



\section{Introduction}

For more than 60 years, the theory of  wave turbulence has been proved to have a wide range applications \cite{Peierls:1993:BRK,Peierls:1960:QTS,hasselmann1962non,hasselmann1974spectral,benney1966nonlinear,kadomtsev1965plasma,zakharov2012kolmogorov,benney1969random}. Among the most  important classes of  wave kinetic equations (WKEs),
3-wave kinetic equations  read
\begin{equation}\label{WeakTurbulenceInitial}
	\begin{aligned}
		\partial_tf(t,k) \ =& \ \mathcal{Q}_{3w}[f](t,k), \\\
		f(0,k) \ =& \ f_0(k),
	\end{aligned}
\end{equation}
in which  $f(t,k)$ is the nonnegative wave density  at  wavenumber $k\in \mathbb{R}^N$, $N \ge 2$; $f_0(k)$ is the initial condition.  The quantity $\mathcal{Q}[f]$ is of the form 
\begin{equation}\label{def-Qf}\mathcal{Q}_{3w}[f](k) \ = \ \iint_{\mathbb{R}^{2N}} \Big[\mathcal R_{k,k_1,k_2}[f] -\mathcal R_{k_1,k,k_2}[f] -\mathcal R_{k_2,k,k_1}[f] \Big] k^Nk_1d^Nk_2, \end{equation}
with $$\begin{aligned}
	\mathcal	R_{k,k_1,k_2} [f]:=  |V_{k,k_1,k_2}|^2\delta(k-k_1-k_2)\delta(\omega -\omega_{1}-\omega_{2})(f_1f_2-ff_1-ff_2), 
\end{aligned}
$$
with the short-hand notation $f = f(t,k)$, $\omega = \omega(p)$ and $f_j = f(t,k_j),$ $\omega_j = \omega(k_j)$, for wavenumbers $k$, $k_j$, $j\in\{1,2\}$. The function $\omega(k)$ is the dispersion relation of the wave system. This type of  wave kinetic equations has several applications from ocean waves, acoustic waves, gravity capillary waves to Bose-Einstein condensates and many others (see \cite{hasselmann1962non,hasselmann1974spectral,nazarenko2019wave,zakharov1968stability,zakharov1967weak,Binh1,PomeauBinh,zakharov1965weak} and references therein). In the isotropic case, we identify  $f(t,k)$ with $f(t,\omega)$ and  the  isotropic 3-wave kinetic equation now has the following form
\begin{align}\label{EE1Colm}
	\begin{split}
		&	\partial_t f(t,\omega) \ = \ {Q}[f](t,\omega), \ \ \ \ \omega\in\mathbb{R}_+,\\
		&	f(0,p) \ = \ f_0(p),\\
		&	{{Q}}[f](t,\omega) \ = \ \int_0^\infty\int_0^\infty \big[R(\omega, \omega_1, \omega_2)-R(\omega_1,\omega, \omega_2)-R(\omega_2, \omega_1, \omega) \big]{d}\omega_1{d}\omega_2, \\\
		& R(\omega, \omega_1, \omega_2):=  \delta (\omega-\omega_1-\omega_2)
		\left[ U(\omega_1,\omega_2)f_1f_2-U(\omega,\omega_1)ff_1-U(\omega,\omega_2)ff_2\right]\,,
	\end{split}
\end{align}
where $U$ satisfies $| U(\omega_1,\omega_2) |  \ = \ (\omega_1\omega_2)^{\gamma/2},$ in which $\gamma$ is a non-negative constant which plays an important role in the sequel.

In \cite{soffer2019energy}, the authors show that if we define the energy  of the solution  \eqref{EE1Colm} as $
g(t,\omega) \ = \ \omega f(t,\omega)$, the energy is conserved. It has been proved in \cite{soffer2019energy} that $g$ can be decomposed into two parts
\begin{equation}
	\label{Decomposition} 
	g(t,\omega) \ = \  \bar{g}(t,\omega) \ + \ \tilde{g}(t)\delta_{\{\omega=\infty\}},
\end{equation}
where $\bar{g}(t,\omega)\ge 0$ is a non-negative function, and $\tilde{g}(t)\delta_{\{\omega=\infty\}}$, is  a measure.   At time $t=0$, we assume $\bar{g}(0,\omega)={g}(0,\omega)$ and $\tilde{g}(0)=0$. It has been proved in \cite{soffer2019energy}, that there exists {\it infinitely many blow-up times} \begin{equation}
	\label{Decomposition1} 0< t^*_1<t^*_2<\cdots<t_n^*<\cdots,\end{equation}
such that \begin{equation}
	\label{Decomposition2}  \bar{g}(t_1^*,\omega)>\bar{g}(t_2^*,\omega)>\cdots>\bar{g}(t_n^*,\omega)> \cdots \to 0,\end{equation}
and  \begin{equation}
	\label{Decomposition3} 0<\tilde{g}(t_1^*)<\tilde{g}(t_2^*)<\cdots<\tilde{g}(t_n^*)< \cdots\end{equation}
In the limit  $t\to\infty$, all of the energy will be accumulated to the  measure $\tilde{g}(t)\delta_{\omega=\infty}$, while the function  will vanish $\tilde{g}(t)\delta_{\omega=\infty}\to0$. We refer to time $t^*_1$ as the {first blow-up} time.
We define $\chi_{[0,R]}(\omega)$ be a cut-off function of $\omega$ on the finite domain $[0,R]$, the {\it multiple blow-up time } phenomenon \eqref{Decomposition}-\eqref{Decomposition1}-\eqref{Decomposition2}-\eqref{Decomposition3}, with the decay rate $\mathcal{O}(\frac{1}{\sqrt{t}})$, can be observed as the decay of the total energy on any finite interval $[0,R]$
\begin{equation}
	\label{Decomposition4}
	\int_{0}^R g(t,\omega)d\omega\ = \	\int_{\mathbb{R}_+}\chi_{[0,R]}(\omega)   g(t,\omega)d\omega\le \  \mathcal{O}\Big(\frac{1}{\sqrt{t}}\Big) \mbox{   as  } t\to\infty,
\end{equation}
for all truncated parameter $R$.
Inequality \eqref{Decomposition4} simply means that the energy of the solution will move away from any truncated finite interval $[0,R]$ as $t\to\infty$ with the rate $\mathcal{O}\Big(\frac{1}{\sqrt{t}}\Big).$ This theoretical finding has been numerically verified via a Finite Volume Scheme in \cite{waltontranFVS}.

In the  important works \cite{connaughton2009numerical,connaughton2010aggregation,connaughton2010dynamical}, several numerical experiments were designed to  investigate time dependent solutions of  isotropic 3-wave equations.
We refer to \cite{waltontranFVS} for a detailed comparison between the different results of \cite{connaughton2009numerical,connaughton2010aggregation,connaughton2010dynamical,soffer2019energy,waltontranFVS}. Below, we recall a brief comparison. 
The works \cite{connaughton2009numerical,connaughton2010aggregation,connaughton2010dynamical} and \cite{soffer2019energy,waltontranFVS} complement each other as they consider very different scenarios of the solutions of  \eqref{EE1Colm}.  The works \cite{soffer2019energy,waltontranFVS} focuses on the finite capacity case $(\gamma>1)$, under the condition that the energy of the solution $f$ of \eqref{EE1Colm} is conserved in time
\begin{equation}\label{ColmEnergy2a}
	\int_{0}^\infty \omega f(t,\omega)d\omega=\mbox{ constant}.
\end{equation}
and shows that   there exist an infinite series of blow up times \eqref{Decomposition}-\eqref{Decomposition1}-\eqref{Decomposition2}- \eqref{Decomposition3} (or  inequality \eqref{Decomposition4}). The works \cite{connaughton2009numerical,connaughton2010aggregation,connaughton2010dynamical} focus on both cases - finite capacity   ($\gamma>1$) and    infinite capacity  $(0\le \gamma\le 1)$. In the finite capacity case ($\gamma>1$),  the solution is studied before the first blow-up time $t<t_1^*$, rigorously proved later in  \cite{soffer2019energy} under the self-similar assumption  
\begin{equation}\label{Ansartz}
	f(t,\omega)\approx s(t)^a F\left(\frac{\omega}{s(t)}\right).
\end{equation} and it is also  assumed that the energy  grows linearly in time 
\begin{equation}\label{ColmEnergy2}
	\int_{0}^\infty \omega f(t,\omega)d\omega=Jt.
\end{equation}
Solving for the  self-similar profiles before the first blow-up times $t_1^*$ is an interesting direction of research. We mention, for instance,   \cite{bell2017self}, where a self-similar profile of the solution for a  different finite capacity system  - the Alfven wave
turbulence kinetic equation - is computed  before the first blow-up time $t_1^*$. In addition, a recently published paper \cite{semisalov2021numerical} presents a numerical method for  for a collision integral 4-wave kinetic equation based on Chebyshev approximations.    

Concerning the analysis of 3-wave kinetic equation, we mention the work \cite{RumpfSofferTran}, where  a 3-wave kinetic equation, derived from the elastic beam wave equation on the lattice, has been studied.  The global existence of 3-wave kinetic equations has been investigated in \cite{GambaSmithBinh,nguyen2017quantum} and the link  to  reaction networks has also been pointed out in \cite{CraciunBinh,CraciunSmithBoldyrevBinh}. 

{\it The goal of our current paper is to present an alternative numerical method, based on deep learning of the solution of the conservative form of equation \eqref{EE1Colm} introduced in \cite{waltontranFVS}. }
It is difficult to exaggerate the impact and scope of machine learning within applied mathematics and scientific computing more broadly.  A very active branch of scientific machine learning is the development of neural network approximations of the solutions to partial differential equations.  As opposed to data-driven discovery of the dynamics of a physical process in which known governing equations may not be used, a physics-informed neural network (PINN) \cite{raissi,lulu} trains the neural network by minimizing the residual of the governing PDE along with minimization terms for the initial and boundary conditions.  For completeness, we describe briefly the general idea.

Let $\mathcal{L}$ be a linear or non-linear (integro-)differential operator.  Let $(t,x)\in[0,T]\times\Omega\subset \mathbb{R}^+\times \mathbb{R}^d$ for $d\geq 1$ and $\Gamma = \partial \Omega$ smooth.  Assume for simplicity that $u(t,x)\in L^2((0,T);\Omega)$ is the solution to the evolution equation
\begin{equation}\label{eqn::evol_eqn}
	\begin{aligned}
		(\partial_t + \mathcal{L})u = 0,\\
		\mathcal{B}u = u_\Gamma, \\
		u(0,x) = u_0(x),
	\end{aligned}
\end{equation}
where $\mathcal{B}$ denotes a boundary operator and $u_0\in L^2(\Omega)$.

If we wish to approximate the solution to \eqref{eqn::evol_eqn} by a neural network $n(t,x;\theta)$ where $\theta\in \Theta$ the set of weights and biases of the neural network, then a common approach is to first define the residual operator, $\mathcal{R}$, such that $$n(t,x;\theta)\in \mathcal{H}_{\mathcal{R}}=\{\varphi(t,x) \; |\; \varphi(t,x) \in L^2((0,T);\Omega), \; \mathcal{R}\varphi(t,x) \in L^2((0,T);\mathcal{R})\},$$
with $$ \| \varphi(t,x) \|_{\mathcal{H}_{\mathcal{R}}} = \| \varphi(t,x) \|_{L^2((0,T);\Omega)} + \| \mathcal{R}\varphi(t,x) \|_{L^2((0,T);\Omega)}, $$
for 
\begin{equation}
	\mathcal{R} = (\partial_t + \mathcal{L}). 
\end{equation}
We then define the functional 
\begin{equation}\label{eqn:J}
	J_\theta[n](t,x) = \| \mathcal{R}n(t,x;\theta)\|^2_{L^2((0,T);\Omega)}+\|\mathcal{B}n(t,x;\theta)-u_\Gamma\|^2_{L^2(\Gamma)} + \|n(0,x;\theta)-u_0 \|^2_{L^2(\Omega)}.
\end{equation}

To obtain the neural network approximation to $u(t,x)$, one then solves the following stochastic optimization problem
\begin{equation}\label{eqn:min_prob}
	\theta^* = \arg\min_{\theta}J_\theta[n](t,x),
\end{equation}
for the (probably local and non-unique) minimizer $\theta^*$ from which we obtain $$n(t,x;\theta^*)\approx u(t,x).$$

In solving \eqref{eqn:min_prob}, one has many choices even when one excludes the usual hyper-parameter tuning that comes with selecting an architecture.  In the absence of empirical measurements, a major choice to be made is how one selects inputs $(t,x)\in(0,T)\times\Omega$ for the residual and initial and boundary conditions.  A common approach is to take samples from the uniform distribution and approximate the integrals in \eqref{eqn:J} by taking expectations (\cite{raissi, lulu} and many, many others).  As in \cite{dgm}, one may use knowledge about the underlying distribution of the sample space to select a sampling distribution which more closely mimics empirical measurements and may be more relevant to the approximated solution.  

In the present article, we choose neither of the previous approaches and opt for a quasi-Monte Carlo (qMC) approximation of the integrals in \eqref{eqn:J} selecting a low-discrepency sequence as our sample points, the Sobol sequence \cite{sobol} (see figure \ref{fig:test2_batch_samples}).

We find the qMC approach to be appropriate for WKEs for two reasons.  The first is that given the highly non-local nature of the collision term, the approximation of this operator is expensive and when coupled with the number of samples necessary to obtain an accurate approximation of the collision term via MC methods, we find a uniform sampling procedure to be impractical for our needs.  The latter method introduces interesting challenges in that the distribution with which we would like to sample the wavenumber domain in order to take expectations is in fact the solution we want to find in the first place.  

Thus a qMC approach which requires fewer sample points and being deterministic requires no knowledge of the underlying distribution is a simple and natural choice.  For a more detailed analysis of qMC methods in neural network approximations to pdes with applications to fluid dynamics the reader is referred to the works \cite{mishra1, mishra2}. 

\subsubsection*{Outline}
We now give a short overview of the remainder of the article.  To build confidence in our approach, we present a neural network (NN) approximation to the solution of a Smoluchowski coagulation equation (SCE) in section 2.  We choose this examples because, 1) as has been discussed elsewhere (\cite{connaughton2010aggregation, soffer2019energy}) the SCE can be considered as a special case of a 3-WKE and 2) unlike for the 3-WKE the SCE has a known analytic solution which we may compare our approximation against. 
In section 3, we present results for a neural network representation of the non-stationary solution to a 3-WKE.  The results are compared with previously derived theoretical results for the decay rate of the total energy in any finite interval of the wavenumber domain.  As a means of validation for the neural network model, we give solutions for the same 3-WKE obtained via a finite volume scheme (FVS).

\section{The Smoluchowski Coagulation Equation}\label{sec:smol}
\par
Before giving results for wave kinetic equations, we provide a check of the method on a similar type of equation, the Smoluchowski coagulation equation, which has a known analytic solution in contrast to wave kinetic equations.  Comparisons of the SCE and 3-WKE can be found in \cite{soffer2019energy, connaughton2010aggregation, waltontranFVS} and references therein.
\par
The SCE may be written in the form \cite{Fil04}
\begin{equation}\label{smol}
	\begin{aligned}
		v\partial_t f(t,v) = -\partial_v \mathcal{Q}_{Smol}[f](t,v) && (t,v) \in \mathbb{R}^+\times\mathbb{R}^+,\\
		f(0,v) = f_0(v)\geq 0 && v\in\mathbb{R}^+:=(0,\infty),
	\end{aligned}
\end{equation}
where the $f(t,v)\geq 0$ gives the density of particles at time $t$ with volume $v$ and 
\begin{equation}\label{smol_collision_operator}
	\begin{aligned}
		\mathcal{Q}_{Smol}[f](t,v) = \int^v_0\int^\infty_{v-v_1}a(v_1,v_2)v_1f(t,v_1)f(t,v_2)\mathrm{d}v_2\mathrm{d}v_1,
	\end{aligned}
\end{equation}
where the kernel is given by $a(v_1,v_2)=v_1v_2$ in what follows.  Let us define
\[
m(t,v) = vf(t,v),
\]
then the total volume can be defined as 
\[
\mathcal{V}(t) = \int^\infty_0 m(t,v) dv.
\]
We will consider the case where 
\begin{equation}\label{smol_IC}
	f_0(v) = \frac{e^{-v}}{v},
\end{equation}
i.e.,
\begin{equation}
	m_0(v) = e^{-v},
\end{equation}
and the analytic solution corresponding to this initial condition is given by \cite{Fil04}
\begin{equation}\label{smol_analytic}
	f(t,v) = e^{-Tv}\frac{I_1(2vt^{1/2})}{v^2t^{1/2}},
\end{equation}
where
\[
I_1(v) = \frac{1}{\pi}\int^\pi_0e^{v\cos\theta}\cos\theta\mathrm{d}\theta,
\]
is the modified Bessel function of the first kind and 
\[
T = 
\begin{cases}
	1 + t & t\leq T_{gel}\\
	2t^{1/2} & t > T_{gel}
\end{cases},
\]
in equation \eqref{smol_analytic} above.
\par 
This SCE undergoes a gelation phenomenon at time $T_{gel} = 1$, that is,the particle number density is entirely concentrated at $v=\infty$ (again, see \cite{Fil04} and the references therein). This also means that $\mathcal{V}(t)=1$ if $t\in[0,T_{gel})$ and $\mathcal{V}(t) = t^{-1/2}$ if $t\geq T_{gel}$. 
\par 

\subsection{NN Representation of the SCE}
To begin, instead of \eqref{smol}, we solve the evolution equation for $m(t,v)$.  Then,
\begin{equation}\label{smol_m}
	\begin{aligned}
		\partial_t m(t,v) = -\partial_v \mathcal{Q}_{Smol}[m](t,v) && (t,v) \in \mathbb{R}^+\times\mathbb{R}^+,\\
		m(0,v) = m_0(v) = v f_0(v)\geq 0 && v\in\mathbb{R}^+,
	\end{aligned}
\end{equation}
and the Smoluchowski collision operator is then written as 
\begin{equation}\label{smol_collision_operator_m}
	\begin{aligned}
		\mathcal{Q}_{Smol}[m](t,v) = \int^v_0\int^\infty_{v-v_1}a(v_1,v_2)m(t,v_1)\frac{m(t,v_2)}{v_2}\mathrm{d}v_2\mathrm{d}v_1.
	\end{aligned}
\end{equation}

Let $m(t,v;\theta)$ be the neural network approximation to the solution of \eqref{smol_m} with $\theta$ denoting the weights and biases of the neural network.  The NN approximation is computed for $(t,v)\in[0,T]\times[0,R]$ for some time $T>0$ and truncation value $R>0$. 
\\
To compute the NN approximation $m(t,v;\theta)$, the problem \eqref{smol_m} is recast as the following functional minimization problem \cite{bensoussan1982lectures,bensoussan2022machine}
\begin{equation}\label{smol_min}
	m(t,v;\theta^*) = \min_{\theta\in\Theta} J_{Smol}[m](t,v;\theta),
\end{equation}
with $\theta^*\in \Theta$ a minimizing set of parameters and $J_{Smol}[m](t,v;\theta)$ is 
\begin{equation}\label{smol_loss}
	J_{Smol}[m](t,v;\theta) = \| \mathcal{R}m(t,v;\theta) \|^2_{L^2((0,T]\times[0,R])} + \| m(0,v;\theta) - m_0(v) \|^2_{L^2([0,R])}, 
\end{equation}
where $\mathcal{R}$ is the residual operator defined by
\begin{equation}\label{smol_res}
	\mathcal{R}m(t,v;\theta) = \partial_t m(t,v;\theta) +\partial_v \mathcal{Q}[m](t,v;\theta).
\end{equation}
The functional $J_{Smol}[m]$ and collision term $\mathcal{Q}_{Smol}[m]$ are discretized via a Quasi-Monte Carlo method (\cite{mishra1}, \cite{mishra2})).  Sample points for the functional are drawn from the Sobol sequence (\cite{caflisch}) in the unit square which are then mapped to the discrete set $S \subset (0,T]\times[0,R]$.  For the collision term, for each $v$, we draw samples from the Sobol sequence in the unit square as before,  but map the variables $(v_1,v_2)$ to the discrete intervals $V_1 \subset [0,v]$ and, for each $v_1$, $V_2 \subset [v-v_1,R]$, respectively.    \\
The discretized collision operator is then given by
\begin{equation}\label{smol_collision_m_discrete}
	\mathcal{\hat{Q}}_{Smol}[m](t,v;\theta) = \frac{v}{|V_1|}\sum_{v_1\in V_1}v_1 m(t,v_1;\theta) \Bigg(\frac{|R+v_1-v|}{|V_2|}\sum_{v_2\in V_2} m(t,v_2;\theta)\Bigg),
\end{equation}
which leads to the semi-discretized residual
\begin{equation}\label{smol_res_disc}
	\hat{r}(t,v;\theta) = \partial_t m(t,v;\theta) +\partial_v \mathcal{\hat{Q}}_{Smol}[m](t,v;\theta),
\end{equation}
and so the semi-discrete functional $\hat{J}_{Smol}[m](t,v;\theta)$ to be minimized is
\begin{equation}\label{smol_loss_discrete}
	\hat{J}_{Smol}[m](t,v;\theta) = \frac{1}{|S|}\sum_{(t,v)\in S} \hat{r}(t,v;\theta)^2 + \frac{1}{|S_0|}\sum_{v\in S_0} (m(0,v;\
	\theta) - m_0(v))^2,
\end{equation}
where the set $S_0 \subset S$, $|S_0|\leq|S|$ denotes a subset of the volume samples in $S$.

\section{A 3-Wave Kinetic Equation}
In \cite{waltontranFVS}, a new identity for the energy of the solutions to 3-WKEs \cite{soffer2019energy} is presented which takes the form of a conservation law \cite{Fil04}.  This is the equation we shall study in the present paper due to the simplicity provided by this form of the equation.  Specifically, there is no need to compute the resonant manifolds of the system, though we emphasize that the method does not depend on this simplification.  

Let us write the equation to be solved. The equation is equivalent to \eqref{EE1Colm} and is explained in \cite{soffer2019energy,waltontranFVS}
\begin{equation}\label{wke}
	\begin{aligned}
		\partial_t g(t,p) = p\partial_p\mathcal{Q}[g](t,p) && (t,p)\in \mathbb{R}^+\times\mathbb{R}^+,\\
		g(0,p) = g_0(p) \geq 0 && (t,p)\in \{0\}\times\mathbb{R}^+,
	\end{aligned}
\end{equation}
where the collision term is given by
\begin{equation}\label{wke_coll}
	\begin{aligned}
		\mathcal{Q}[g](t,p) = -2\int^p_0\int^p_0 (p_1,p_2)^{\frac{\gamma}{2}-1}g_1g_2\chi_{p,1,2}\mathrm dp_{2,1} + \int^\infty_0\int^\infty_0(p_1,p_2)^{\frac{\gamma}{2}-1}g_1g_2\chi_{p,1,2}\mathrm dp_{2,1}, 
	\end{aligned}
\end{equation}
where we have used the notation $g_i = g(t,p_i)$ for $i=1,2$, $\mathrm{d}p_{2,1} = \mathrm{d}p_2\mathrm{d}p_1$ and $\chi_{p,1,2} = \chi\{p < p_1 + p_2 \}$ where $\chi\{ A \}$ is the set characteristic function of some set $A$.  The parameter $\gamma$ is the degree of the kernel as discussed in the introduction and $\gamma =2 $ 
in the present work which corresponds to acoustic wave systems.

As in the previous section, we define a neural network approximation, $g(t,p;\theta)$ to be a solution to the optimization problem \cite{bensoussan1992stochastic}
\begin{equation}\label{wke_min}
	g(t,p;\theta^*) = \min_{\theta\in\Theta} J[g](t,p;\theta),
\end{equation}
with $\theta^*\in \Theta$ a minimizing set of parameters of the functional 
\begin{equation}\label{wke_loss}
	J[g](t,p;\theta) = \| r(t,p;\theta) \|^2_{L^2(\mathbb{R}^+\times\mathbb{R}^+)} + \| g(0,p;\theta)-g_0(p)\|^2_{L^2(\mathbb{R}^+)},
\end{equation}
where, again, $\mathcal{R}$ denotes the residual operator of the evolution equation defined by
\begin{equation}\label{wke_res}
	\mathcal{R}g(t,p;\theta) = \partial_t g(t,p;\theta) - p\partial_p\mathcal{Q}[g](t,p;\theta).
\end{equation}

The functional \eqref{wke_loss} is approximated via a Quasi-Monte Carlo method with sample points drawn from the Sobol sequence in the unit square and then transformed to some truncated rectangle of the time, wavenumber domain, i.e. we generate the set $W \sim Sobol([0,T]\times[0,R])$, for $T,R >0$ truncation parameters of the time and wavenumber domain, respectively.  The residual \eqref{wke_res} is approximated similarly where the set of sample points for $(p_1,p_2)$ are given by $P_1 \sim Sobol(0,p)$ for each $p$, and $P_2\sim Sobol(p-p_1,p)$ for each $p,p_1$ in the first term of \eqref{wke_coll} and $\hat{P}_2 \sim Sobol(p-p_1,R)$ for each $(p,p_1)$.  

To see how these sample sets are defined, note that the collision term can be rewritten as 
\begin{equation}
	\begin{aligned}
		\mathcal{Q}[g](t,p) = -2\int^p_0\int^p_{p-p_1}(p_1p_2)^{\frac{\gamma}{2}-1}g_1g_2\mathrm{d}p_{2,1} + \int^p_0\int^R_{p-p_1}(p_1p_2)^{\frac{\gamma}{2}-1}g_1g_2\mathrm{d}p_{2,1},
	\end{aligned}
\end{equation}
where we have applied the truncation parameter of the wavenumber domain to the second expression and enforced the restrictions $p-p_1 < p_2$ in both terms to satisfy the characteristic set function and $p_1 < p$ to guarantee positivity of the variable $p_2$ in the second term.  Thus, using the above collision operator and sample sets defined in the previous paragraph, we can define the discrete collision operator to be
\begin{equation}
		\mathcal{\hat{Q}}[g](t,p) = \frac{-2p}{|P_1||P_2|}\sum_{p_1\in P_1}p_1^{\frac{\gamma}{2}}g_1\Bigg(\sum_{p_2\in P_2}p_2^{\frac{\gamma}{2}-1}g_2 \Bigg) + \frac{p}{|P_1||\hat{P}_2|}\sum_{p_1\in P_1}p_1^{\frac{\gamma}{2}-1}g_1\Bigg[(R-p+p_1)\Big(\sum_{p_2\in \hat{P}_2}p_2^{\frac{\gamma}{2}-1}g_2 \Big)\Bigg].
\end{equation}
%

Using the discrete collision operator, we can define the semi-discrete residual by
\begin{equation}\label{wke_res_discrete}
	\hat{\mathcal{R}}g(t,p;\theta) = \partial_t g(t,p;\theta) - p\partial_p \mathcal{\hat{Q}}[g](t,p;\theta),
\end{equation}
and therefore we can minimize the semi-discrete functional 
\begin{equation}\label{wke_loss_discrete}
	\hat{J}[g](t,p);\theta) = \frac{1}{|W|}\sum_{(t,p)\in W}\hat{\mathcal{R}}g(t,p;\theta)^2 + \frac{1}{|W_0|}\sum_{p\in W_0}(g(0,p;\theta)-g_0(p))^2,
\end{equation}
where $W_0\subset W$, $|W_0|\leq |W|$ contains only sample points for the wavenumber.

\section{Numerical Results}
In this section, we provide numerical results for the SCE with initial condition \eqref{smol_IC} and compare with the anaylitc solution \eqref{smol_analytic}.   A 3-WKE with initial condition \eqref{wke_IC} is then solved and the decay rate of the total energy is compared with the results obtained in \cite{soffer2019energy}.  The same initial condition \eqref{wke_IC} is then used to obtain a solution  with the finite volume scheme developed in \cite{waltontranFVS}.  The results of the two methods are discussed. 

All deep learning tests were implemented in TensorFlow \cite{tensorflow} with Keras \cite{keras}.  Training was accelerated using multiple GPUs utilizing Horovod \cite{horovod}.  All computations were performed on SMU's computing cluster MII.  The finite volume computations were performed in Matlab.

\subsection{Test 1}

A simple, yet effective, architecture is chosen.  Namely, it is enough to employ a feedforward network with only two hidden layers with 128 hidden units in each layer utilizing sigmoidal activation functions \cite{Hornik}.  Given its simplicity, we can write the neural network (NN) out explicitly as 
\begin{equation*}
	m(t,v;\theta) = \mathbf{W}\sigma(\mathbf{W}_2\sigma(\mathbf{W}_1\vec{x}+\mathbf{b}_1) + \mathbf{b}_2) + b,
\end{equation*}
with the entries of $\mathbf{W}_i,\mathbf{b}_i \sim \mathcal{N}(0,1)$ for $i=1,2$, $\mathbf{W}, b \sim \mathcal{N}(0,1)$ and $\sigma(\cdot)$ denoting the sigmoid function.
To fully discretize the functional \eqref{smol_loss_discrete}, we employ TensorFlow's built-in automatic differentiation method.

\begin{figure}
	\begin{center}
		\includegraphics[width = 13cm, height = 11cm]{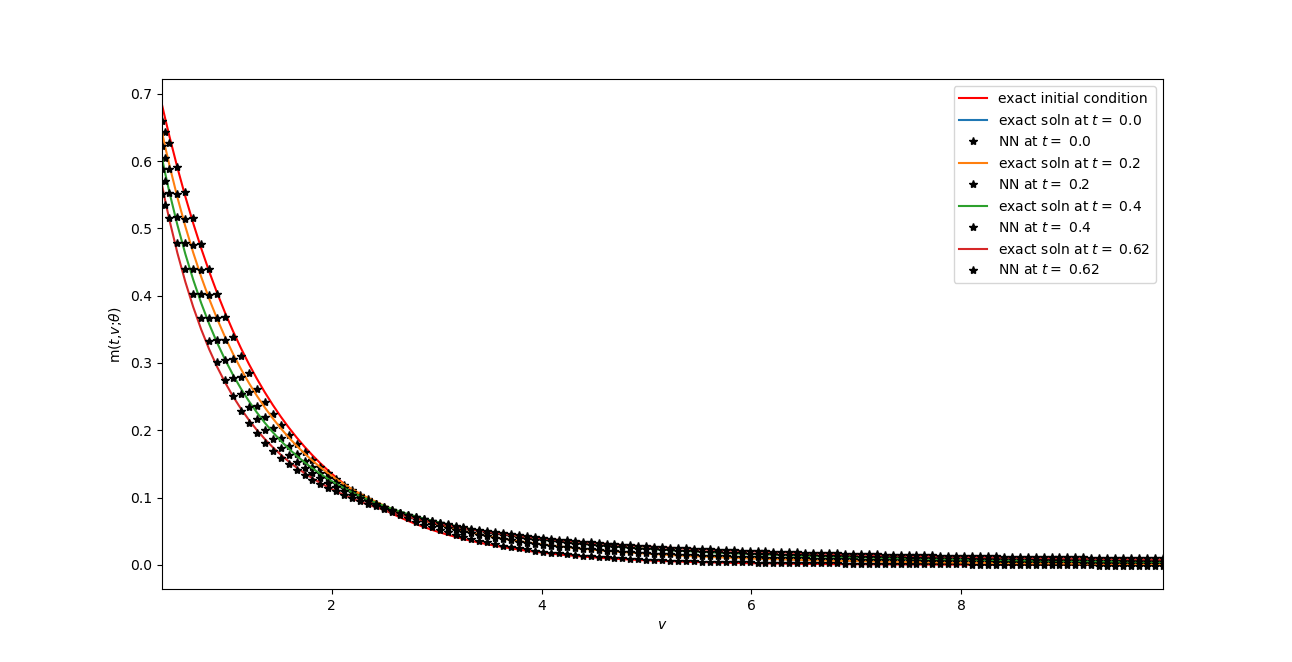}
		\caption{\small Neural Network (NN) approximation to the Smoluchoswki equation \eqref{smol_m}.  The solid lines denote the analytic solution \eqref{smol_analytic} and $*$ markers denote the NN approximation at equivalent snapshots in time. }
		\label{fig::smol_net_solution}
	\end{center}
\end{figure}

The training samples were drawn from the rectangle $(0,T]\times[0,R]$ using the strategy described above for $T=0.8$ and $R=8$.  To achieve the results shown in Figure \ref{fig::smol_net_solution}, only the first 32 Sobol points were used to minimize the residual term (16 time samples and 16 volume samples) and only the first 16 Sobol points were needed to train on the initial data.

Figure \ref{fig::smol_net_solution} shows the neural network predictions on unseen data within the training interval.  The inputs to the neural network to produce the plot were the first 128 Sobol points transformed to the volume training interval and $t=0.0, 0.2, 0.4$ and $0.62$.  We see the NN solution gives a good approximation to the analytic solution.  

\begin{figure}
	\begin{center}
		\includegraphics[scale = 0.8]{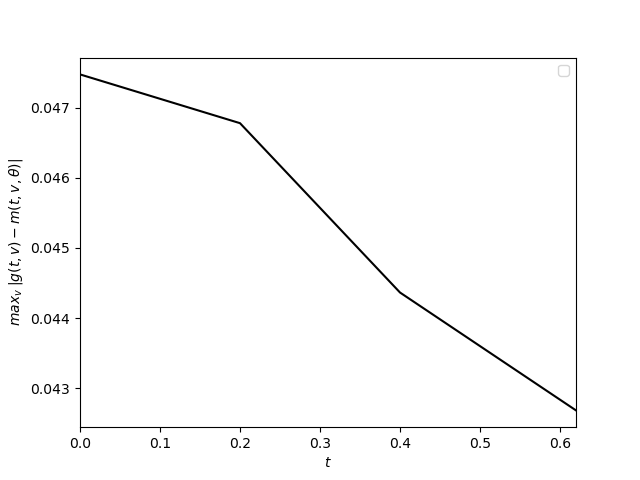}
		\caption{\small Sup norm of the error at $t = 0.0, 0.2, 0.4$ and $ 0.62$ where the maximum is taken over the volume domain $v\in [0,10]$.  This shows the accuracy of the method on unseen data points within the training interval. }
		\label{fig::smol_sup_error}
	\end{center}
\end{figure}

\begin{figure}
	\begin{center}
		\includegraphics[scale = 0.8]{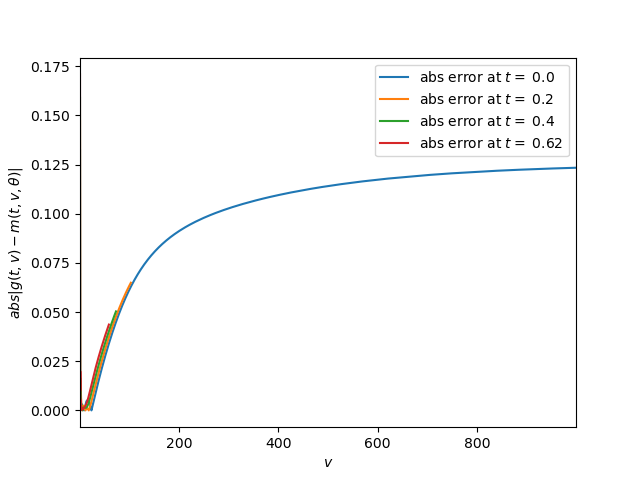}
		\caption{\small Absolute error of the NN approximation versus the analytic solution over the volume interval $v\in[0,1000]$ with $2^{13}$ volume samples given as input to the NN.  We see a slow growth in the error for large volume numbers. }
		\label{fig::smol_abs_error}
	\end{center}
\end{figure}

In Figure \ref{fig::smol_sup_error}, we quantify the match of the NN approximation to the analytic solution by giving the sup norm of the error at each snapshot. 
Here, we supply the first $2^{10}$ Sobol points within the volume training interval and use the same points in time $t = 0.0, 0.2, 0.4, 0.62$.  We note here that while the NN performs well on unseen points within the volume interval $v\in[0,R]$, for $v> R$ we generalization error is larger.  This is to be expected given the small training sample size used (16 volume points).  We show the absolute error in Figure \ref{fig::smol_abs_error} for each snapshot in time over the interval $v\in[0, 1000]$ with $2^{13}$ Sobol points in this volume interval. This result informed our decision to apply a batch sampling procedure for the wave kinetic equation, the details of which are discussed below. 

\subsection{Test 2}

We choose the initial condition
\begin{equation}\label{wke_IC}
	g_0(p) = \sqrt{\frac{7}{2\pi}}e^{-\frac{7(p-2)^2}{2}} .
\end{equation}
The neural network was trained using samples from the Sobol sequence in the rectangle $W\sim Sobol((0,T]\times[0,R])$ for $T=10$ and $R=10$. 
\begin{figure}
	
	\includegraphics[scale=0.45]{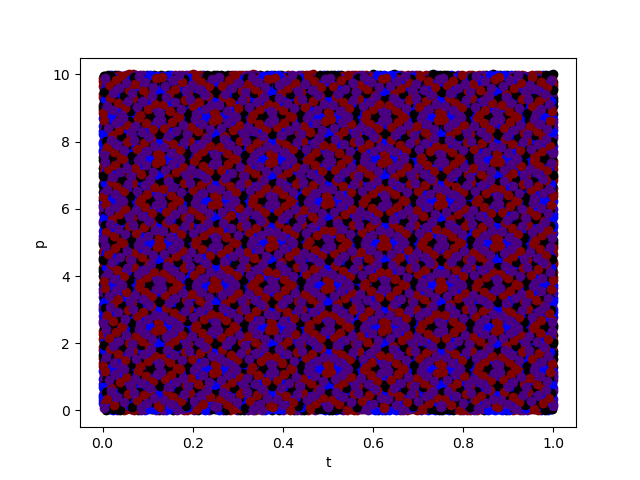}    \includegraphics[scale = 0.45]{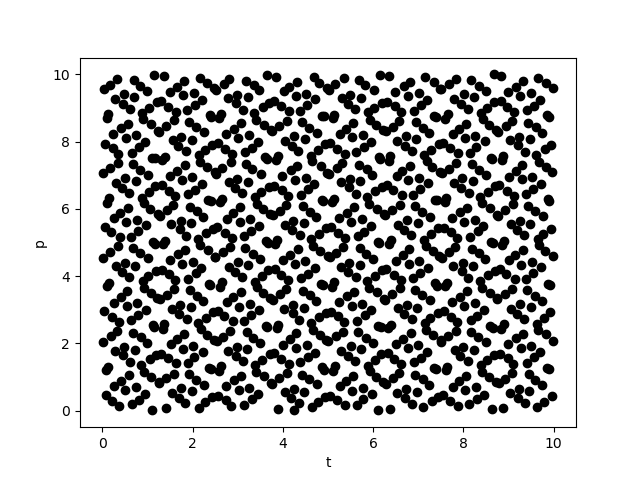}\\  \includegraphics[scale =0.45]{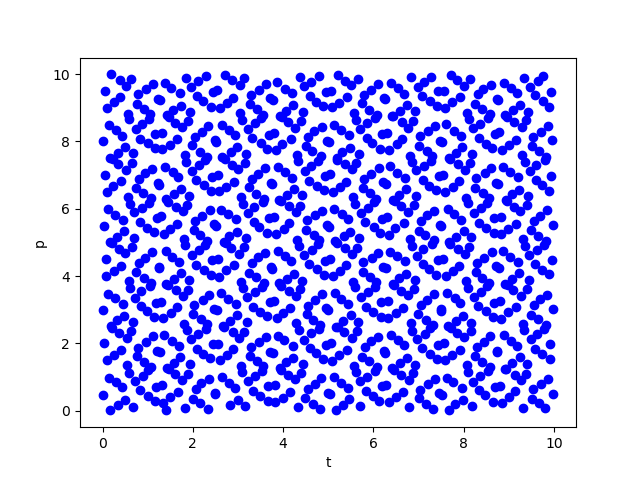} \includegraphics[scale = 0.45]{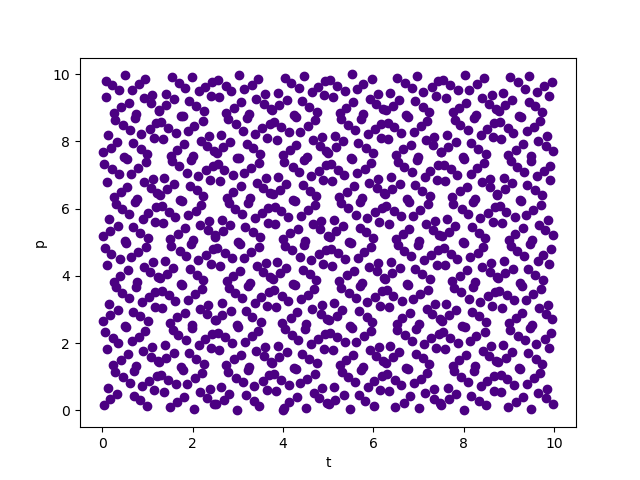}
	\caption{\small Batched samples of Sobol points. (Top Left)  All sample points. (Top Right and Bottom  Row) Example batches of sampled points.}
	\label{fig:test2_batch_samples}
	
\end{figure}
As illustrated in figure \ref{fig:test2_batch_samples}, the size of the sample set $|W|$ was $2^{15}$, which was broken up into smaller training batches.

The architecture was again chosen to have 2 hidden layers, each with 128 units and sigmoidal activation functions.  The loss was again minimized using tensorflows implementation of ADAM.  The collision terms were approximated with the sampling procedure described above with each sample comprised of 32 Sobol points.  We consider only the case $\gamma = 2$ here.

\begin{figure}
	\centering
	\includegraphics[scale=0.76]{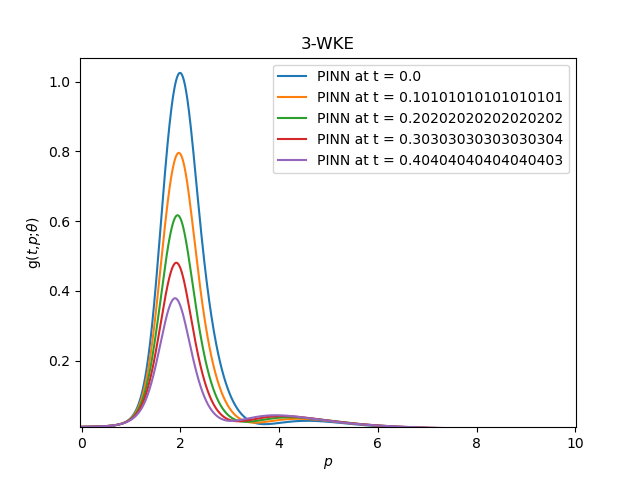}
	\includegraphics[scale = 0.31]{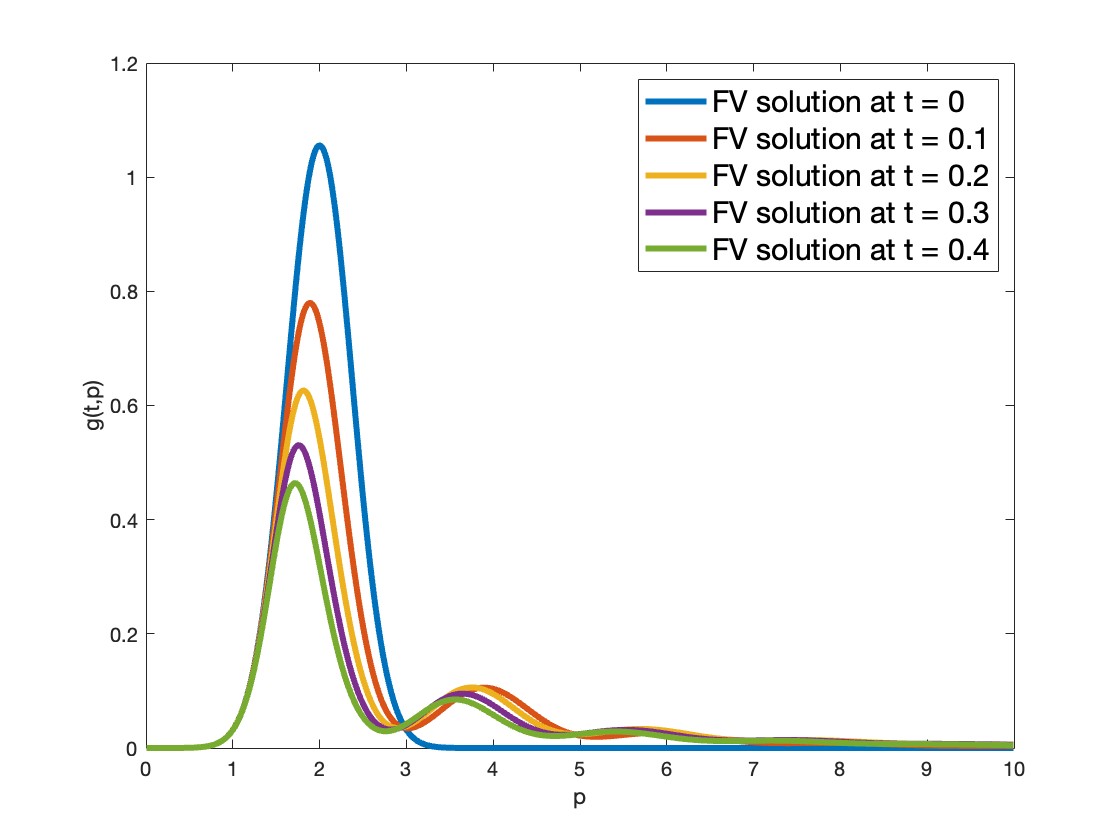}
	\caption{Top Picture: A few snapshots of the NN approximation corresponding to initial condition \eqref{wke_IC}. Bottom Picture: Comparative snapshots of the FVS solution for the same initial condition \eqref{wke_IC}.  }
	\label{fig:test2_snapshots}
\end{figure}

We show a few early snapshots of the approximated solution in Figure \ref{fig:test2_snapshots}.  The behaviour seems consistent with predictions in that the $L^\infty$ norm decays as time increases.  This is more clearly evidenced by the total energy discussed in the following paragraph. In Figure \ref{fig:test2_snapshots}, we also present a comparison with the solution computed the Finite Volume Scheme of \cite{waltontranFVS}. This comparison will be described later in this section, around equation \eqref{FVSST}-\eqref{FVSscheme}.

\begin{figure}
	\centering
	\includegraphics[scale = 0.765]{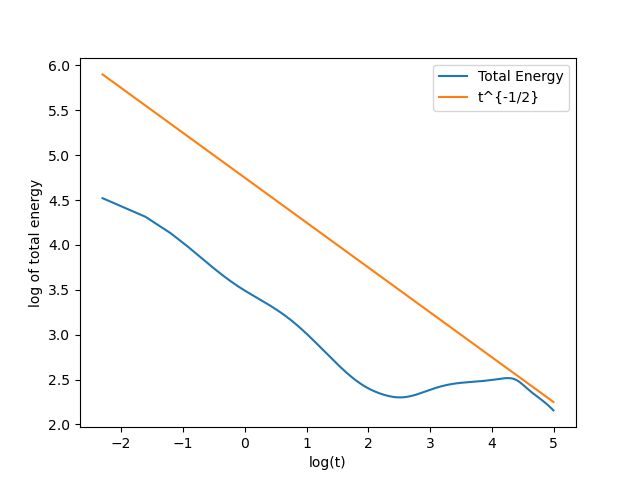}
	\includegraphics[scale = 0.31]{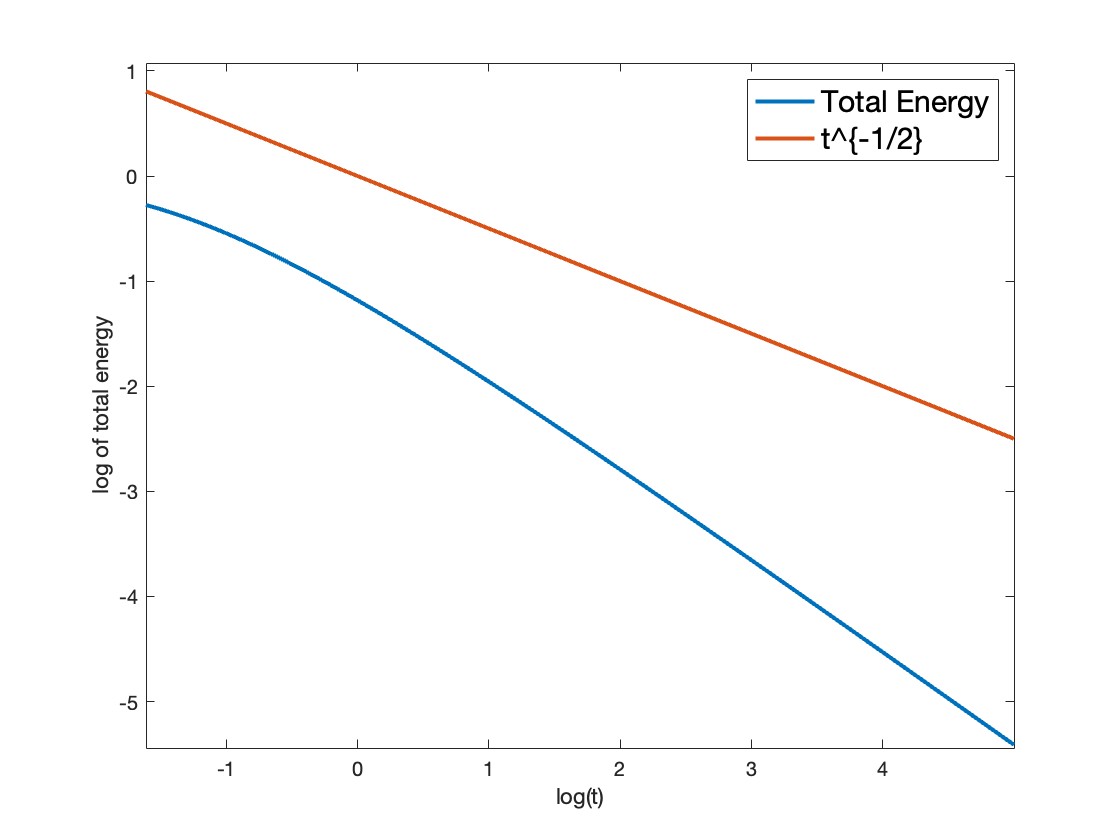}
	\caption{Top Picture: Log-Log plot of the total energy corresponding to initial condition \eqref{wke_IC} as predicted up to $t=148$ the  NN.  Bottom Picture: The total energy computed for $t\in[0,148]$ with the FVS. }
	\label{fig:test2_total_energy}
\end{figure}

The total energy of the solutions was computed and the log of the total energy is plotted against the log of $t$. The comparison between the decay of the numerical solution is made with the theoretical decay rate of $t^{-1/2}$, obtained in \cite{soffer2019energy} (see \eqref{Decomposition4}), in Figure \ref{fig:test2_total_energy}.  The total energy was predicted up to $t=148$.  From Figure \ref{fig:test2_total_energy}, the energy  decays and the decay is in good agreement with the theoretical rate obtained in \cite{soffer2019energy}.

\begin{figure}
	\centering
	\includegraphics[scale = 0.65]{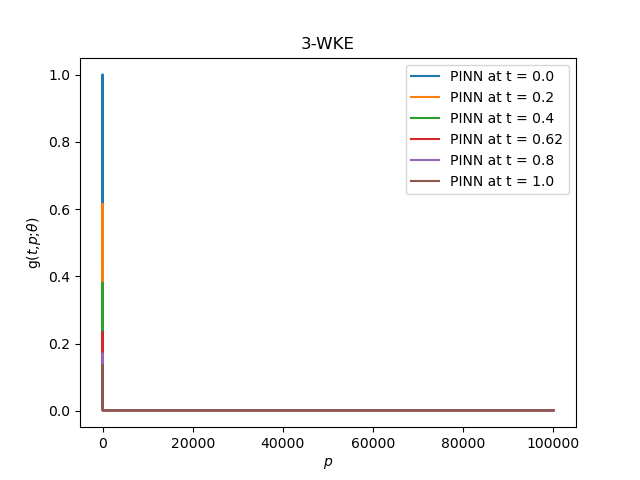}\\
	\includegraphics[scale = 0.23]{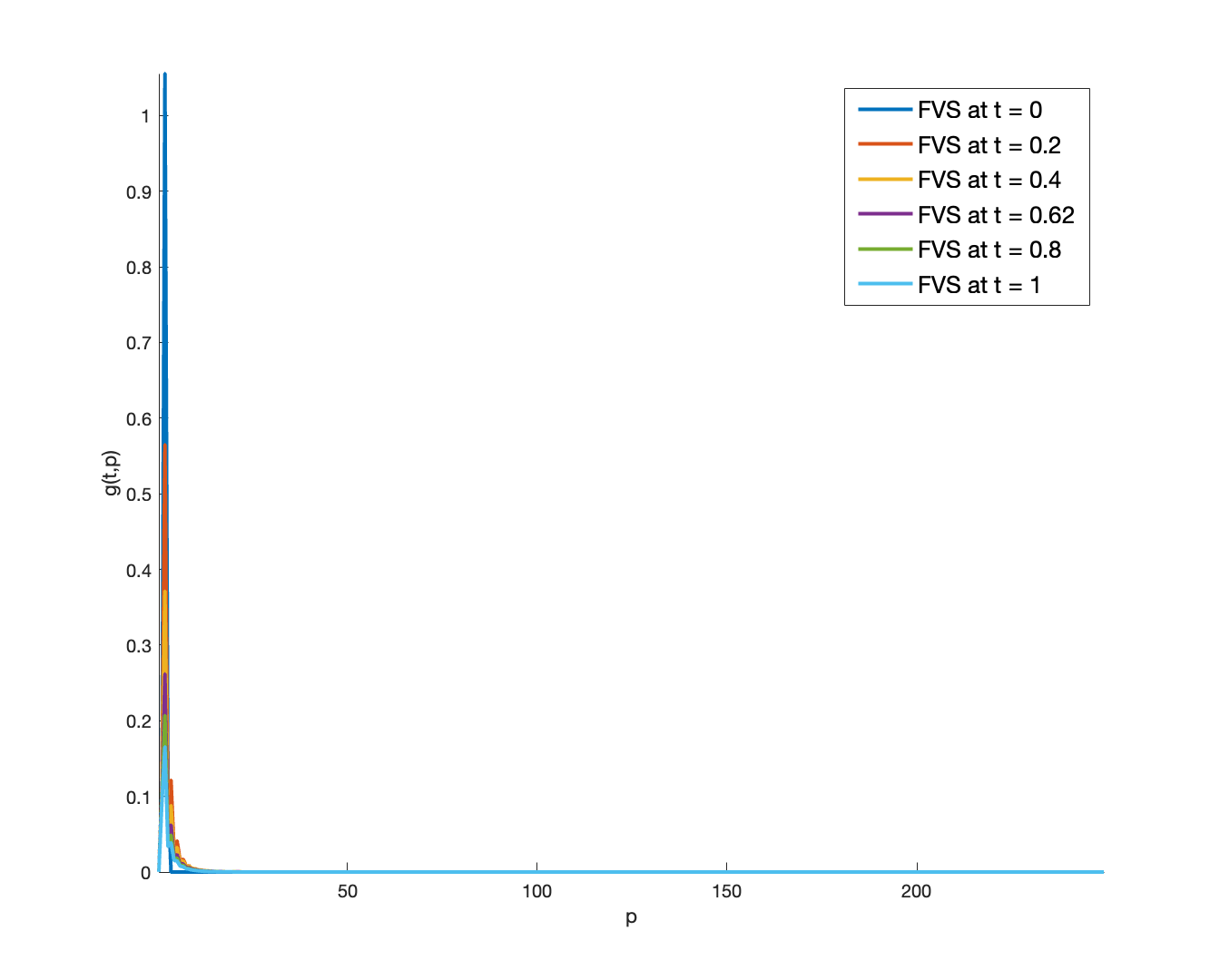}
	\caption{Top picture: Neural network prediction corresponding to initial condition \eqref{wke_IC} for wavenumbers up to 1e+06. Bottom picture: The finite volume solution for wavenumbers up to $p=250$, the largest value for which stability is maintained without further decreasing $\Delta t$.  The two figures highlight a key strength of the presented method, in that the neural network approximation is consistent with more traditional solvers but is able to produce results for computationally prohibitive values of the wavenumber for traditional methods while maintaining positivity and stability.    }
	\label{fig:test2_p1m}
\end{figure}

The neural network is able to make reasonable predictions for very large, unseen in training, samples in the wave number domain.  In Figure \ref{fig:test2_p1m}, we see the predicted solution for wavenumber values up to 1e+06.  In contrast, the most we could push for reasonable predicitons in time was $t=148$ as previously reported in the results shown in figure \ref{fig:test2_total_energy}.  For better predictions in time, we needed to enforce more dense sampling in the time domain.  Thus, the neural network was not only trained on the set $W$ for T= 10, but also for $T=5$ and $T=2$ with the same number of sample points, $2^{15}$.  

\subsection*{Comparison with a Finite Volume Scheme}

As mentioned earlier, to validate the results, we provide a comparison with another model.  Here we use the finite volume scheme presented in \cite{waltontranFVS}, which we briefly outline below.

The discretization of the wavenumber domain $p\in [0,R]$, is given as follows. Let $i\in \{0,1, 2, \ldots , M \} = I^M_h$, with $h \in (0,1)$ fixed and $M=M(h)$.  Define the set  $S_d = \{0,\ldots, R\}$  to be the discretization of the interval $[0,R]$.  Let
\begin{equation}\label{FVSST}
	S_d = \{p_{i+1/2}\}_{i\in I^M_h},\;
	\{p_i\}_{i\in I^M_h\setminus\{0\}} = \frac{p_{i+1/2} +  p_{i-1/2}}{2},\; \{\Delta p_i\}_{i\in I^M_h\setminus\{0\}} = p_{i+1/2} - p_{i-1/2} \leq h,
\end{equation}
define the faces, pivots and step-size respectively, with $p_{1/2} = 0$ and $p_{M+1/2} = R$. For simplicity, we use a uniform grid and set $h=0.01$ which leads to
\[
\begin{aligned}
	S_d = \{ih\}_{i\in I^M_h}, && \{p_i\}_{i\in I^M_h\setminus\{0\}} = \frac{h}{2}(2i-1), && \{\Delta p_i\}_{i\in I^M_h\setminus\{0\}} = h \in (0,1).
\end{aligned}
\]
The set $T_N = \{0,\ldots,T\} $ with $N+1$ nodes and $T=148$ is the maximum time.  We fix the time step to be $\Delta t = \frac{T}{N} = 0.005$, and denote by $t_n = \Delta t\cdot n$ for $n\ \in \{0,\ldots,N\}$.  We approximate equation \eqref{wke} by 
\begin{equation}\label{FVSscheme}
	g^{n+1}(p_{i}) = g^n(p_{i})  +\lambda_i\Big( Q^n_{i+1/2}\Big[\frac{g}{p} \Big] - Q^n_{i-1/2}\Big[\frac{g}{p} \Big] \Big),
\end{equation}
where $\lambda_i = \frac{p_{i}\Delta t }{\Delta p_i}$, and 
\[
Q^n_{i+1/2}\Big[\frac{g}{p} \Big] - Q^n_{i-1/2}\Big[\frac{g}{p} \Big] = -2\Big(Q^n_{1,i+1/2}\Big[\frac{g}{p} \Big]-Q^n_{1,i-1/2}\Big[\frac{g}{p} \Big]\Big) + \Big(Q^n_{2,i+1/2}\Big[\frac{g}{p} \Big]-Q^n_{2,i-1/2}\Big[\frac{g}{p} \Big]\Big),
\]
with
\begin{equation}\label{J}
	Q^n_{1,i+1/2}\Big[\frac{g}{p} \Big] = \sum^{i}_{m=1}\Delta p_m \frac{g^n(p_m)}{p_m} \Bigg(\sum^i_{j=1}\Delta p_j \frac{g^n(p_j)}{p_j}a(p_m, p_j)\chi\Big\{p_{i+1/2} < p_m +p_j \Big\}
	\Bigg),
\end{equation}

\begin{equation}\label{C}
	Q^n_{2,i+1/2}\Big[\frac{g}{p} \Big] =  \sum^M_{m=1}\Delta p_m \frac{g^n(p_m)}{p_m} \Bigg(\sum^M_{j=1}\Delta p_j \frac{g^n(p_j)}{p_j}a(p_m, p_j)\chi\Big\{p_{i+1/2} < p_m+p_j\Big\}\Bigg),
\end{equation}
where we have used the midpoint rule to approximate the integrals in equation \eqref{wke} and we choose an explicit time stepping method. 

The initial condition \eqref{wke_IC} is approximated by
\[
g^0(p_{i}) = \frac{1}{\Delta p_i}\int^{p_{i+1/2}}_{p_{i-1/2}}g_0(p) \mathrm{d}p \approx g_0(p_i),
\]
by again employing the midpoint rule.

We draw the reader's attention again to Figure \ref{fig:test2_snapshots} where a few comparative snapshots of the solution are provided. Qualitatively, the two models seem to agree and capture the main features of the theorized behavior of solutions.  Namely, an evacuation of the energy within any finite interval.  Further, as is typical of simple feedforward architectures, the finer oscillations seen in the FVS solutions appear to be averaged out in the NN solution as expected.  For our purposes, the dynamics captured by the simple feed-forward architecture we have employed here are enough to confirm the theory presented in \cite{soffer2019energy}.  We leave it to a future work to investigate more complicated architectures.  For example, the architecture described in \cite{lizuo} is able to capture highly oscillatory solutions for the stationary Navier-Stokes equations. 

In Figure \ref{fig:test2_total_energy} (bottom) we see the decay rate of the total energy as provided by the computed solution of the FVS, which is compared with the {\em predicted} values by the NN (top). Both models are in good agreement with the theorized bound on the rate of decay \cite{waltontranFVS, soffer2019energy}, though the FVS appears to capture a slightly faster rate of decay.    

Looking back to Figure \ref{fig:test2_p1m}, we make a comparison with the computed solution of the finite volume method for large wavenumber values with predictions provided by the neural network.  Here, we have increased $h$ to $0.8$ and $R=250$ while keeping the timestep fixed at $\Delta t = 0.005$ as in the previous figures. The two figures highlight a key strength of the presented deep learning method, in that the neural network approximation is consistent with more traditional solvers but is able to produce results for computationally prohibitive values of the wavenumber for traditional methods, while maintaining positivity and stability.  Indeed, constructing positivity preserving schemes for PDEs is a very important direction of research \cite{hu2013positivity,huang2019positivity,huang2019third}. 
 For the FVS \eqref{FVSscheme}, the positivity and thus stability is lost for wavenumber values larger than 250 with the timestep set at $\Delta t = 0.005$.  In order to preserve the positivity of the solutions produced by the FVS \eqref{FVSscheme}, the CFL condition is restrictive and the time step $\Delta t$ needs to be chosen sufficiently small as shown in Proposition 3.1 of  \cite{waltontranFVS}. As thus, being able to preserve the positivity of the solutions is indeed a very important feature of the presented deep learning approximation.


\section{Conclusions}
\label{sec:conclusions}

We present a deep learning approximation, stochastic optimization based,  method for the 3-wave kinetic equation, studied theoretically in \cite{soffer2019energy} and numerically in \cite{waltontranFVS}. We first apply the method to a Smoluchowski coagulation equation with multiplicative kernel for which an analytic solution exists. The deep learning method is proved to give a good approximation of the analytic solution.  Next, the learning approach is then used to approximate the solution of the 3-wave kinetic equation. The deep learning approximation is tested and proved to be as good as the Finite Volume approximation introduced in \cite{waltontranFVS} for the same equation and both approximations are in good agreement with the theoretical results of  \cite{soffer2019energy}.

\section*{Acknowledgments}
Computational resources for this research were provided by SMU’s Center for Research Computing. S.W. would like to thank Thom Hagstrom for many useful discussions.



\bibliographystyle{elsarticle-num} 
\bibliography{ml_wke.bib}





\end{document}